\documentclass[12pt,a4paper]{article}
\usepackage{german}
\usepackage{latexsym}
\usepackage{amsfonts}
\usepackage{amssymb}
\usepackage{amsbsy}
\usepackage{mathrsfs}
\usepackage{amsmath}
\usepackage{epsfig}
\usepackage{theorem}
\usepackage{color}
\input{xy}
\xyoption{all}

\def\colim{\mathop{\rm colim}}
\def\hocolim{\mathop{\rm hocolim}}

\def\Top{\mathcal{T}op}
\def\Sp{\mathcal{S}p}
\def\Ass{\mathcal{A}ss}
\def\sC{\mathcal{C}}

\def\scA{{\mathcal{A}}}

\def\scC{{\mathcal{C}}}

\def\scK{{\mathcal{K}}}

\def\op{\mathop{\rm op}}

\makeatletter
\@addtoreset{equation}{section}

\makeatother

\theoremstyle{change}
\theorembodyfont{\rm}
\newtheorem{prop}{Proposition:}[section]
\newtheorem{theo}[prop]{Theorem:}
\newtheorem{cond}[prop]{Condition:}
\newtheorem{coro}[prop]{Corollary:}
\newtheorem{rema}[prop]{Remark:}
\newtheorem{defi}[prop]{Definition:}
\newtheorem{leer}[prop]{}

\newtheorem{lem}[prop]{Lemma:}
\newenvironment{proof}{\par\noindent\textbf{Proof:}}{\hfill\ensuremath{\Box}}

\begin{document}
\title{\textbf{On the Multiplicative Structure of Topological Hochschild Homology}}
\author{\textbf{M. Brun, Z. Fiedorowicz, R.M.~Vogt}}
\date{\textbf{17.06.2007}}
\maketitle

\paragraph*{Abstract}
We show that the topological Hochschild homology $THH(R)$ of an $E_n$-ring
spectrum $R$ is an $E_{n-1}$-ring spectrum. The proof is based on the fact
that the tensor product of the operad $\Ass$ for monoid structures and the
the little $n$-cubes operad $\mathcal{C}_n$ is an $E_{n+1}$-operad, a result
which is of independent interest.

\vspace{2ex}
In 1993 Deligne asked whether the Hochschild cochain complex of an
associative ring has a canonical action by the singular chains of the
little 2-cubes operad. Affirmative answers for differential graded algebras
in characteristic $0$ have been found by 
Kontsevich and Soibelman \cite{KS}, Tamarkin \cite{Tam1} and
\cite{Tam2}, and Voronov \cite{Vor}. A more general
proof, which also applies to associative ring spectra is due to McClure and
Smith \cite{MS}. In \cite{Kon} 
Kontsevich extended  Deligne's question: Does the Hochschild cochain complex of
an $E_n$ differential graded algebra carry a canonical $E_{n+1}$-structure?

We consider the dual problem: Given a ring $R$ with additional structure,
how much structure does the topological Hochschild homology $THH(R)$ of $R$
inherit from $R$? The close connection of $THH$ with algebraic $K$-theory and
with structural questions in the category of spectra make multiplicative
structures on $THH$ desirable.

In his early work on topological Hochschild homology of functors with smash
product B"okstedt proved that $THH$ of a commutative such functor is a
commutative ring spectrum (unpublished). The discovery of associative,
commutative and unital smash product functors of spectra simplified the
definition of $THH$ and the proof of the corresponding result for
$E_\infty$-ring spectra considerably (e.g. see \cite{MSV}).

In this paper we morally prove

\textbf{Theorem A:}
For $n\geq 2$, if $R$ is an $E_n$-ring spectrum 
then $THH(R)$ is an $E_{n-1}$-ring
spectrum.

The same result has been obtained independently by Basterra and Mandell 
using different techniques \cite{BM}.

Why ``morally''? To define $THH(R)$ we need $R$ to be a strictly
associative spectrum. In general, $E_n$-structures do not have a strictly
associative substructure. So we have to replace $R$ by an equivalent
strictly associative ring spectrum $Y$, whose multiplication extends to an
$E_n$-structure. Then the statement makes sense for $Y$. 
Here is a more  precise
reformulation of Theorem A.
 
\textbf{Theorem B:}
Let $R$ be an $E_n$-ring spectrum. Then there are $E_n$-ring spectra 
$X$ and $Y$ and
maps of $E_n$-ring spectra
$$
\xymatrix{
Y & X \ar[l]\ar[r] & R
}
$$
which are homotopy equivalences of underlying spectra such that
the $E_n$-structure on $Y$ extends a strictly associative ring
structure and the topological Hochschild homology 
$THH(Y)$ inherits an $E_{n-1}$-ring structure from $Y$.

Theorem B is an easy consequence of the universality of the
$W$-construction of Boardman-Vogt \cite{BV2}, \cite{Vogt}
and
an interchange result involving the operad structures
of the operad $\scA ss$ codifying monoid structures and the little $n$-cubes operad
$\scC _n$. The interchange is codified by the tensor product of operads (for terminology
see \ref{1_6}). Our key result will be

\textbf{Theorem C:} $\scA ss\otimes \scC_n$ is an $E_{n+1}$-operad satisfying Condition
\ref{1_2} below.

In \cite{F} this theorem has been announced and the main ideas for a proof have been
sketched. 
Here we will include a detailed proof by depicting the spaces 
$(\scA ss\otimes \scC_n)(k)$ as iterated
colimits of diagrams of contractible spaces over posets. The diagrams of this iterated
colimit combine to give a diagram over a Grothendieck construction, whose realization
will turn out to be an $E_{n+1}$-operad.

Since this iterated colimit construction might be of use in other
cases we give a formal definition in an appendix.

In this paper we
will be working in the categories $\Top$, $\Top^\ast$, and $\Sp$ of $k$-spaces,
based $k$-spaces, and $S$-module spectra in the sense of \cite{EKMM}. Hence
spectrum will always mean $S$-module spectrum.

\section{$\mathbf{E_n}$-operads}
In this section we recall the basic definitions underlying the statement of Theorem C.

\begin{defi}\label{1_1}
An \textit{operad} 
is a topologically enriched strict symmetric mo\-noidal
category $(\mathcal{B},\oplus,0)$ such that 
\begin{itemize}
\item $ob\:\mathcal{B}=\mathbb{N}$ and $m\oplus n=m+n$
\item
$\mathcal{B}(m,n)=\coprod\limits_{r_1+\ldots +r_n=m}\mathcal{B}(r_1,1)\times
\ldots\times
\mathcal{B}(r_n,1)\times_{\Sigma_{r_1}\times\ldots\times
\Sigma_{r_n}}\Sigma_m $
\end{itemize}
\end{defi}

Since the morphism spaces $\mathcal{B}(m,n)$ are determined by the spaces
$\mathcal{B}(r,1)$ we usually write $\mathcal{B}(r)$ for $\mathcal{B}(r,1)$
and work with them. A \textit{map of operads} $f:\mathcal{B}\to\mathcal{C}$
is a continuous strict symmetric monoidal functor such that $f(n)=n$ for
all $n\in\mathbb{N}$. It is called a \textit{weak equivalence} if
$f:\mathcal{B}(n)\to\mathcal{C}(n)$ is a homotopy equivalence of spaces for
all $n$. We call it a \textit{$\Sigma$-equivalence} if these map are
$\Sigma_n$-equivariant homotopy equivalences.\\
An operad $\mathcal{B}$ is called $\Sigma$-\textit{free} if
$\mathcal{B}(r)\to\mathcal{B}(r)/\Sigma_r$ is a numerable principal
$\Sigma_r$-bundle for all $r$.

For technical reasons we sometimes require

\begin{cond}\label{1_2}
$\{id\}\subset \mathcal{B}(1)$ is a closed cofibration.
\end{cond}

\begin{defi}\label{1_3}
Let $\mathcal{C}_n$ denote the little $n$-cubes operad \cite[Chap. 2, Expl.
5]{BV1}. An $E_n$-\textit{operad} is an operad $\mathcal{B}$
for which there exists a sequence of $\Sigma$-equivalences of operads
$$
\xymatrix{
\mathcal{B}=\mathcal{B}_0 \ar[r]^(0.6){f_0}& \mathcal{B}_1 & \ldots
  \ar[l]_(0.4){f_1} & \ar[r] & \mathcal{B}_r & \mathcal{C}_n\ar[l]_(0.4){f_r}
}
$$
\end{defi}

\begin{leer}\label{1_4}
Let $\mathcal{O}pr$ denote the category of operads. In \cite[Chap.
III]{BV2} Boardman-Vogt constructed as continuous functor
$$
\xymatrix{
W: \mathcal{O}pr \ar[r] & \mathcal{O}pr
}
$$
together with a natural transformation
$$
\xymatrix{
\varepsilon: W \ar[r] & Id
}
$$
taking values in $\Sigma$-equivalences.
It can
be interpreted as a cofibrant replacement of $\mathcal{B}$ \cite{Vogt}. 
In particular, given a diagram of maps of operads
$$
\xymatrix{
 &  \mathcal{C}\ar[d]^g\\
W\mathcal{B} \ar[r]^f &  \mathcal{D}
}
$$
such that $\mathcal{B}$ satisfies Condition \ref{1_2} and
$g$ is a $\Sigma$-equivalence,
then there exists a lift $h:W\mathcal{B}\to \mathcal{C}$ up to homotopy
through operad maps, and $h$ is unique up to homotopy through operad maps.
If $\mathcal{B}$ is also $\Sigma$-free, the same holds if $g$ is only a weak
equivalence \cite[3.17]{BV2}.
\end{leer}

This implies

\begin{prop}\label{1_5}
Let $\mathcal{B}$ and $\mathcal{C}$ be $E_n$-operads such that
$\mathcal{B}$ satisfies Condition (\ref{1_2}). Let 
$$
\xymatrix{
\mathcal{B}=\mathcal{B}_0 \ar[r]^(0.6){f_0}& \mathcal{B}_1 
& \ldots \ar[l]_(0.4){f_1} & \ar[r] & \mathcal{B}_r 
& \mathcal{C}\ar[l]_(0.4){f_r}
}
$$
be a sequence of weak equivalences connecting them. Then there is a diagram
of weak equivalences
$$
\xymatrix{
&& && W\mathcal{B}
\ar[ddllll]_{\varepsilon_\mathcal{B}} \ar[ddll]_{g_1}
\ar[ddrrrr]^{g_{r+1}} \ar[ddrr]^{g_r}
&& &&
\\
\\
\mathcal{B} \ar[rr]^{f_0} && \mathcal{B}_1
&& \ldots \ar[rr] \ar[ll]_{f_1} &&
\mathcal{B}_r && \mathcal{C}\ar[ll]_{f_r}
} 
$$
commuting up to homotopy through operad maps. The $g_i$ are unique up to
homotopy through operad maps.
In particular, there exists a weak equivalence
$W\mathcal{C}_n\to\mathcal{C}$.

$E_n$-structures are closely related to $n$ interchanging
$(E_1=A_\infty)$-structures. Let $\mathcal{C}$ and $\mathcal{D}$ be two
operads and $X$ be an object having a $\mathcal{C}$- and a
$\mathcal{D}$-structure. These structured are said to \textit{interchange}
if for each $c\in \mathcal{C}(n)$ the operation
$$
\xymatrix{
c:X^n \ar[r] & X
}
$$
is a $\mathcal{D}$-homomorphism, or equivalently, for each
$d\in\mathcal{D}(m)$ the operation
$$
\xymatrix{
d:X^m \ar[r] & X
}
$$
is a $\mathcal{C}$-homomorphism, i.e. the diagram
\end{prop}

\begin{leer}\label{1_6}
$$
\xymatrix@C=0pt{
(X^m)^n\ar[d]^{d^n} &\cong & (X^n)^m \ar[rrrrrr]^(0.6){c^m} 
&&&&&& X^m \ar[d]^d
\\
X^n\ar@{-}[rrrrrrrr]^(0.5)c  &&&&&&&& X
}
$$ 
commutes for all $c\in\mathcal{C}(n)$ and all $d\in\mathcal{D}(m)$. 
\end{leer}

The resulting structure on $X$ is codified by an operad
$\mathcal{C}\otimes\mathcal{D}$ called the \textit{tensor product} of
$\mathcal{C}$ and $\mathcal{D}$. Formally, $\mathcal{C}\otimes\mathcal{D}$
is the quotient of the categorical sum $\mathcal{C}\oplus\mathcal{D}$ in
$\mathcal{O}pr$ by factoring out relation (\ref{1_6}). For more
details see \cite[p. 40ff]{BV2}.

Theorem C will be proved in Section 4. In \ref{4_12} we will give an explicit chain of 
$\Sigma$-equivalences of operads connecting $(\scA ss\otimes \scC_n)$ with
$\scC_{n+1}$.

\section{Algebraic structures on spectra}
\begin{leer}\label{2_1}
The category $\Sp$ of spectra is enriched over $\Top^\ast$ and complete and
cocomplete in the enriched sense (for details see \cite[Chap. VII]{EKMM}).
If $K$, $L$ are based spaces and $M,N$ are spectra we have a natural
isomorphism
$$
M\wedge (K\wedge L)\cong (M\wedge K)\wedge L
$$
and natural homeomorphisms
$$
\Sp(M\wedge K,N)\cong \Top^\ast(K,\Sp(M,N)\cong \Sp(M,F(K,N))
$$
where $F(K,N)$ is the function spectrum. In particular, 
$$
\xymatrix{
-\wedge K: \Sp \ar[r] &\Sp
}
$$
preserves colimits.
\end{leer}

\begin{leer}\label{2_2}
We can form the based topological \textit{endomorphism operad} $\mathcal{E}nd_M$,
given by
$$
\mathcal{E}nd_M(n) =\Sp(M^{\wedge n}, M)
$$
with the $0$-map as base point, where $M^{\wedge 0}=S$ is the sphere spectrum. 

If $\mathcal{C}$ is any operad in $\Top$ , a $\mathcal{C}$-structure on $M$ is an
operad map
$$
\xymatrix{
\mathcal{C}_+ \ar[r] & \mathcal{E}nd_M
}
$$
where $\mathcal{C}_+(n)=\mathcal{C}(n)_+=\mathcal{C}(n)\sqcup\{\ast\}$ with
basepoint $\ast$. This transforms the topological operad $\mathcal{C}$
into a based topological operad $\mathcal{C}_+$; the monoidal structure
in $\Top^\ast$ is given by the smash product.
Passing to adjoints a $\mathcal{C}$-structure on $M$ is
given by a sequence of maps 
$$
\xymatrix{
\mathcal{C}(n)_+\wedge_{\Sigma_n} M^{\wedge n} \ar[r] & M, \quad
n\in\mathbb{N}
}
$$
satisfying certain conditions due to the fact that
$\mathcal{C}_+\to\mathcal{E}nd_M$ is a symmetric monoidal functor.

$M$ together with a given $\mathcal{C}$-structure is called a
$\mathcal{C}$-\textit{algebra} or $\mathcal{C}$-\textit{ring spectrum}.

To make sense of the interchange diagram (\ref{1_6}) we have to give
$M^{\wedge n}$ a $\mathcal{C}$-structure: If $M$ and $N$ are
$\mathcal{C}$-algebras, then the canonical $\mathcal{C}$-algebra structure
on $M\wedge N$ is given by the maps 

\begin{multline*}
\mathcal{C}(n)_+ \to 
(\mathcal{C}(n)\times \mathcal{C}(n))_+
=\mathcal{C}(n)_+\wedge \mathcal{C}(n)_+\to 
\\
\mathcal{E}nd_M(n)\wedge\mathcal{E}nd_N(n) \to 
\mathcal{E}nd_{M\wedge N}(n)
\end{multline*}

where the first map is induced by the diagonal and the last by the smash
product.
\end{leer}

Finally we will need

\begin{prop}\label{2_3}
If $M_\ast$ is a simplicial $\mathcal{C}$-algebra, then the realization
$|M_\ast|$ inherits a $\mathcal{C}$-algebra structure.
\end{prop}

This follows from the fact that $-\wedge \mathcal{C}(n)_+$ preserves
colimits. For details see \cite[X. 1.3, X. 1.4]{EKMM}.

\section{$\mathbf{THH}$ of $\mathbf{E_n}$-ring spectra}
\begin{defi}\label{3_1}
Let $\mathcal{B}$ and $\mathcal{C}$ be $E_n$-operads, let $R$ be a
$\mathcal{B}$-algebra and $M$ be a $\mathcal{C}$-algebra. An
$E_n$-\textit{ring map} from $R$ to $M$ is a pair $(\alpha,f)$ consisting
of an operad map $\alpha:\mathcal{B}\to\mathcal{C}$ and a homomorphism
$f:R\to M$ of $\mathcal{B}$-algebras, where the $\mathcal{B}$-structure on
$M$ is the pulled back $\mathcal{C}$-structure.  
\end{defi}

\begin{leer}\label{3_2}
Let $R\to M$ be an $E_n$-ring map between $E_n$-ring spectra. In view of
(\ref{1_5}) and Theorem C we may assume that it is a homomorphism
of $W(\mathcal{A}ss \otimes\mathcal{C}_{n-1})$-algebras.
\end{leer}

For any operad $\mathcal{B}$ we have the free algebra functor 
${\mathbb{B}}$
from spectra to $\mathcal{B}$-algebras defined by
$$
{\mathbb{B}}(X) = 
\bigvee\limits_{n\ge 0} \mathcal{B}(n)_+\wedge_{\Sigma_n}
X^{\wedge n}.
$$

We now form the monadic bar constructions \cite[Chap. XII]{EKMM} to obtain a
diagram of $E_n$-ring spectra (here $\mathbb{W}(\mathcal{C})$ stands
for the free algebra functor associated with the operad $W\mathcal{C}$)
$$
\def\objectstyle{\scriptstyle}
\xymatrix{
B({\mathbb{A}ss\otimes\mathbb{C}_{n-1}},
{\mathbb{W}(\mathcal{A}ss\otimes\mathcal{C}_{n-1})},R) \ar[d]
& 
B({\mathbb{W}(\mathcal{A}ss \otimes\mathcal{C}_{n-1})},
{\mathbb{W}(\mathcal{A}ss\otimes\mathcal{C}_{n-1})},R) \ar[d]\ar[r]\ar[l]
& R \ar[d]
\\
B({\mathbb{A}ss\otimes\mathbb{C}_{n-1}},
{\mathbb{W}(\mathcal{A}ss\otimes\mathcal{C}_{n-1})},M )
& 
B({\mathbb{W}(\mathcal{A}ss \otimes\mathcal{C}_{n-1})},
{\mathbb{W}(\mathcal{A}ss\otimes\mathcal{C}_{n-1})},M)\ar[r]\ar[l]
& M
}
$$
and $E_n$-ring maps by (\ref{2_3}). Since $\varepsilon:
W(\mathcal{A}ss\otimes\mathcal{C}_{n-1})\to\mathcal{A}ss\otimes
\mathcal{C}_{n-1}$
is a weak $\Sigma$-equivalence the horizontal maps are homotopy
equivalences of spectra \cite[X. 2.4]{EKMM}. Let
$$
\xymatrix{
f:Y_R \ar[r] & Y_M
}
$$
denote the left vertical $(\mathcal{A}ss\otimes\mathcal{C}_{n-1})$-algebra
homomorphism. In particular, $f$ is a homomorphism of strictly associative,
unital ring spectra, so that $Y_M$ is a $Y_R$-bimodule. We can form the
topological Hochschild homology of $Y_R$ with coefficients in $Y_M$:

\begin{leer}\label{3_3}
Let $Q$ be a monoid in $\Sp$ and $N$ a $Q$-bimodule. Then $THH(Q;N)$ is
defined to be the realization of the simplicial spectrum
$$
\xymatrix{
[n] \ar[r] & THH(Q;N)_n = Q^{\wedge n}\wedge N
}
$$
with the well-known Hochschild boundary and degeneracy maps.

The inclusions of the $0$-skeleton defines a natural map 
$$
\xymatrix{
\eta: N \ar[r] & THH(Q; N).
}
$$
In our situation $THH(Y_R;Y_M)_\ast$ is a simplicial
$\mathcal{C}_{n-1}$-algebra by the interchange relation (\ref{1_6}). Hence
$THH(Y_R; Y_M)$ is a $\mathcal{C}_{n-1}$-algebra and 
$$
\xymatrix{
\eta: Y_M \ar[r] & THH(Y_R; Y_M)
}
$$
is a homomorphism of $\mathcal{C}_{n-1}$-algebras. We obtain the following
generalization of Theorem B.
\end{leer}

\begin{theo}\label{3_4}
Let $f: R\to M$ be an $E_n$-ring map between $E_n$-ring spectra. Then there
is a commutative diagram of $E_n$-ring spectra and $E_n$-ring maps
$$
\xymatrix{
Y_R \ar[d]^{f_Y} & X_R \ar[l]\ar[r]\ar[d]^{f_X} & R\ar[d]^{f}
\\
Y_M & X_M \ar[r]\ar[l] & M
}
$$
with the following properties:
\vspace{-2ex}
\begin{enumerate}
\item The horizontal maps are homotopy equivalences of spectra.
\item $Y_R$ and $Y_M$ are $(\mathcal{A}ss\otimes\mathcal{C}_{n-1})$-algebras
and $f_Y$ is an $(\mathcal{A}ss\otimes\mathcal{C}_{n-1})$-algebra
homomorphism.
\end{enumerate}
\vspace{-2ex}
The second property implies that $THH(Y_R;Y_M)$ is a
$\mathcal{C}_{n-1}$-ring spectrum, and the natural map $\eta:Y_M\to
THH(Y_R; Y_M)$ is a $\mathcal{C}_{n-1}$-algebra homomorphism.
\hfill\ensuremath{\Box}
\end{theo}

\section{Proof of Theorem C}
In general analyzing the homotopy type of the tensor product of operads is
an intractable problem. 
Our strategy is to represent $(\Ass\otimes\mathcal{C}_n)(k)$ as the colimit of a diagram
$$
F_k:\mathcal{K}_{n+1}(k)\longrightarrow \Top
$$
of contractible spaces, indexed by the $k$-th space of a  modification of
Berger's complete graphs operad, 
a poset operad defined below, such
that
\begin{itemize}
\item the diagrams are compatible with the operad structures of $\mathcal{K}_{n+1}$ and 
$\Ass\otimes\mathcal{C}_n$
\item the canonical map $\hocolim F_k\to\colim F_k=(\Ass\otimes\mathcal{C}_n)(k)$ 
is a homotopy equivalence
\end{itemize}

Then the collection of the $\hocolim F_k$ forms an operad $\hocolim F$ 
and we have a chain of weak equivalences of operads
$$
\xymatrix{
|\mathcal{K}_{n+1}|=\hocolim_{\mathcal{K}_{n+1}}\ast &
\hocolim F \ar[l] \ar[r] &
\Ass\otimes\mathcal{C}_n
}
$$
Since the topological realization $|\mathcal{K}_{n+1}|$ of $\mathcal{K}_{n+1}$ is a
$\Sigma$-free topological operad, so is $\hocolim F$.
We will show that $(\Ass\otimes\mathcal{C}_n)$ is $\Sigma$-free, hence both weak 
equivalences are $\Sigma$-equivalences.
Moreover, using a corrected version of Berger's argument in
\cite{Be} we will prove that
$|\mathcal{K}_{n+1}|$ is an $E_{n+1}$-operad. Hence
$\Ass\otimes\mathcal{C}_n$ is an 
$E_{n+1}$-operad, too.

\textit{The modified complete graphs operad} $\scK$: 
A {\it coloring\/} of the complete graph on the set of vertices
$\{1,2,3,\dots,k\}$ is an assignment of
colors to each edge of the graph from the countable set of colors
$\{1,2,3,\dots\}$. A {\it monochrome acyclic orientation\/} of a colored
 complete graph on $k$ vertices is an assignment of direction to each
edge of the graph such that no directed cycles of edges of the same color occur.
The poset $\scK(k)$ has as elements pairs
$(\mu,\sigma)$, where $\mu$ is a coloring and $\sigma$ is a monochrome acyclic
orientation of the complete graph on $k$ vertices.
The order relation on
$\scK(k)$ is determined as follows: we say that
$(\mu_1,\sigma_1)\le (\mu_2,\sigma_2)$ if for
any colored oriented edge $a\stackrel{i}{\longrightarrow}b$ in $(\mu_1,\sigma_1)$
the corresponding edge in $(\mu_2,\sigma_2)$ has
orientation and coloring $a\stackrel{j}{\longrightarrow}b$ with $j\ge i$
or $b\stackrel{j}{\longrightarrow}a$ with $j>i$.  The $n$-th filtration
$\scK_{n}(k)$ is the subposet of $\scK(k)$ where the colorings
are restricted to take values in the subset $\{1,2,3,\dots,n\} $.

The action of the symmetric group $\Sigma_k$ on $\scK(k)$ is via permutation
of the vertices. The composition
$$\scK(k)\times\scK(m_1)\times\scK(m_2)
\times\dots\times\scK(m_k)
\longrightarrow\scK_(m_1+m_2+\dots+m_k)$$
assigns to a tuple of orientations and colorings in
$\scK(k)\times\scK(m_1)\times\scK(m_2)
\times\dots\times\scK(m_k)$
the orientation and coloring obtained by subdividing the set of
$m_1+m_2+\dots+m_k$ vertices into $k$ adjacent blocks containing $m_1$, $m_2$, \dots, $m_k$
vertices respectively. The edges connecting vertices within the $i$-th block
are oriented and colored according to the given element
in $\scK(m_i)$.  The edges connecting vertices between blocks $i$ and $j$
are all oriented and colored according to the corresponding edge in the
given element of $\scK(k)$.

Berger's complete graphs operad $\scK^B_n$ is the suboperad of  $\scK_n$ consisting of
those oriented colored graphs which do not have any cycles, i.e..
 polychromatic cycles are also disallowed for elements in $\scK^B_n(k) $
\begin{leer}\label{4_1} \textbf{Analysis of} $\Ass\otimes \sC_n$: 
By \cite[Thm. 5.5]{BV3} the space $(\Ass\otimes \sC_n)(k)$ is the quotient of
$\Sigma_k\times \scC_n(1)^k$ by the relation
$$
(\pi ;c_1,\ldots , c_k)\sim (\rho ;c_1,\ldots , c_k)
$$
iff $\pi^{-1}(i)<\pi^{-1}(j)$ and $\rho^{-1}(i)>\rho^{-1}(j)$
 imply that $(c_i,c_j)\in \sC_n(2)$. The element $(\sigma, c_1,\ldots, c_k)
\in \Sigma_k\times \scC_n(1)^k$ represents the operation
$(x_1,\ldots,x_k)\mapsto
c_{\sigma 1}(x_{\sigma 1})\cdot\ldots\cdot c_{\sigma k}(x_{\sigma k})$,
where $\cdot $ stands for the monoid multiplication.
This also specifies the operad structure.
\end{leer}

\textit{Observation:} Since $(\Ass\otimes \sC_n)(1)=\sC_n(1)$ and
$\sC_n$ satisfies Condition \ref{1_2}, so does $\Ass\otimes \sC_n$.

Since there is only one color for edges of elements in ${\mathcal{K}}_1(k)$,
the operads ${\mathcal{K}}_1$ and ${\mathcal{K}}^B_1$ coincide. In particular,
the elements in ${\mathcal{K}}_1(k)$ do not contain any cycles.
An orientation with no cycles of the complete graph on the set of vertices 
$\{1,2,\ldots ,k\}$ is a total ordering of $\{1,2,\ldots ,k\}$, which in turn
can be identified with a permutation of $\{1,2,\ldots ,k\}$. 
Hence a representative of an element in $(\Ass \otimes\mathcal{C}_n)(k)$ 
can be identified with an oriented graph $\lambda \in {\mathcal{K}}_1(k)$
together with a labelling of the vertices by elements of $\mathcal{C}_n(1)$.  

To take care of the relation (\ref{4_1}) we enlarge the modified complete graphs operad:
we allow complete graphs with {\it partial monochrome acyclic orientations\/} and 
{\it partial colorings\/}. Such graphs $\lambda'$ are obtained from oriented colored
graphs $\lambda \in {\mathcal{K}}(k)$ by choosing a subset $S$ of the set $E(k)$ of
edges of $\lambda$ and forgetting the orientations and colors of all edges in $S$.
The graph $\lambda''$ obtained from $\lambda$ by forgetting the orientations 
and colors of all edges in $E(k)\backslash S$ is called a \textit{complementary
graph of} $\lambda$. Let $\widehat{\scK}(k)$ denote the poset of all pairs
$(\mu,\sigma)$, where $\mu$ is a partial coloring and $\sigma$ is a partial 
monochrome acyclic
orientation of the complete graph on $k$ vertices obtained from some element in
$\scK(k)$. The order relation is defined as follows: 
$(\mu_1,\sigma_1)\le (\mu_2,\sigma_2)$ if every uncolored unoriented edge
in $(\mu_1,\sigma_1)$ is also uncolored unoriented in $(\mu_2,\sigma_2)$, and for
any colored oriented edge $a\stackrel{i}{\longrightarrow}b$ in $(\mu_1,\sigma_1)$
the corresponding edge in $(\mu_2,\sigma_2)$ is either uncolored unoriented or has
orientation and coloring $a\stackrel{j}{\longrightarrow}b$ with $j\ge i$
or $b\stackrel{j}{\longrightarrow}a$ with $j>i$. 

The symmetric group actions and composition in $\widehat{\scK}$, and the $n$-th 
filtration $\widehat{\scK}_n$ are defined as in $\scK$. 
We shall refer to $\widehat{\scK}$ and its filtrations
as the {\it augmented complete graphs operad\/}.

While the topological realization of $\scK_{n}$ is an
$E_n$-operad, this is not true for $\widehat{\scK}_{n}$:
$|\widehat{\scK}_{n}(k)|$ is equivariantly contractible to the $\Sigma_k$
fixed point specified by the complete graph on $\{1,2,\dots,k\}$ with all its
edges unoriented and uncolored.

If we now label the vertices of $\lambda\in \widehat{\scK}_{1}(k)$ by elements
in $\sC_n(1)$ with the extra condition that the pair of labels $(c,c')$ of the
end points of a non-oriented edge is an element of $\sC_n)(2)$, then $\lambda$ with
its vertex labels $(c_1,\ldots,c_k)$ represents the equivalence class 
in $(\Ass\otimes \sC_n)(k)$ of all
$(\lambda';c_1,\ldots,c_k)$, where $\lambda'\in  \scK_{1}(k)$ is an element from
which $\lambda$ can be obtained by forgetting orientations and colors. 
These labelled augmented complete graphs form an operad 
$\widehat{\mathcal{K}}_1\# \mathcal{C}_n$. 
Its composition is induced by the composition in 
$\widehat{\mathcal{K}}_1$ and the following labelling condition: if the $i$-th
vertex of the element in $\widehat{\mathcal{K}}_1(k)$ has the label $a$
we compose the labels of the vertices of the elements in 
$\widehat{\mathcal{K}}_1(m_i)$ from the left with $a$. 
So $(\widehat{\mathcal{K}}_1 \# \mathcal{C}_n)(k)$ is the disjoint
union of all $A(\lambda),\ \lambda \in \widehat{\mathcal{K}}_1(k)$,
where  
 $A(\lambda)\subset\mathcal{C}_n(1)^k$ denotes the space of possible vertex labels of
 $\lambda$. We obtain that $(\Ass \otimes\mathcal{C}_n)(k)$ is a quotient of
$(\widehat{\mathcal{K}}_1\# \mathcal{C}_n)(k)$. More precisely, the analysis
of $\Ass \otimes\mathcal{C}_n$ of \ref{4_1} can be restated as 

\begin{lem}\label{4_2}
$(\Ass \otimes\mathcal{C}_n)(k)$ is the colimit of the diagram
$$
A:\widehat{\mathcal{K}}_1(k)^{op} \longrightarrow \Top, \quad \lambda\mapsto A(\lambda)
$$
where we consider each poset as a category with a morphism $\lambda_1\to\lambda_2$ whenever 
$\lambda_1\le\lambda_2$.
\hfill$\square$
\end{lem}

Our next step is to depict $A(\lambda)$ as a colimit of contractible 
subspaces. Here $\mathcal{K}_{n+1}$
comes into the picture. We embed $\widehat{\mathcal{K}}_1$ into
$\widehat{\mathcal{K}}_{n+1}$ by changing the color $1$ of the colored edges
of its graphs to $n+1$,
and we define 
$T_k(\lambda)$ be the subposet of $\widehat{\mathcal{K}}_n(k)$ of all $\lambda'$
which are complementary graphs of $\lambda$. More explicitly: let $S\subset E(k)$
be the set of colored edges of $\lambda$; then an element 
$\lambda' \in \widehat{\mathcal{K}}_n(k)$ lies in $T_k(\lambda)$ if there   
is an element $\overline{\lambda}\in \mathcal{K}_{n+1}(k)$ whose edges in $S$ are
oriented as in $\lambda$ and colored $n+1$ and $\lambda'$ is obtained from
$\overline{\lambda}$ by forgetting orientations and colors of the edges in $S$.

We define $n$ strict order relations on $\mathcal{C}_n(1)$ as follows:
Let $c_1,c_2\in\mathcal{C}_n(1)$ and let $(x_1,\ldots,x_n)$ be the highest corner of
$c_1$ and $(y_1
,\ldots,
y_n)$ the lowest corner of $c_2$. For  $1\le i\le n$ we define
$$
c_i <_i c_2 \quad \textrm{ iff }\quad x_i\le y_i
$$

For each $\mu\in\widehat{\mathcal{K}}_n(k)$ we define a closed subspace
$H(\mu)\subset \mathcal{C}_n(1)^k$
by
$$
H(\mu)=\{(c_1,\ldots,c_k)\in \mathcal{C}_n(1)^k; \; c_p<_ic_q \textrm{ if } p\stackrel{i}{\to} q
\textrm{ in }\mu\},
$$
and we have a functor
$$
F_k(\lambda): T_k(\lambda) \longrightarrow\Top,\qquad \lambda'\mapsto
\bigcup\{ H(\mu);\ \mu\in T_k(\lambda)\ \mu\le\lambda' \}
$$
where the union is taken in $\mathcal{C}_n(1)^k$.


\begin{lem}\label{4_3} $A(\lambda) = \colim F_k(\lambda)$. Moreover, for
any element $\alpha \in T_k(\lambda)$ the restriction of this colimit to the subposet
$\mathcal{P}=\{\beta\in T_k(\lambda); \ \beta<\alpha \}$ is a subspace of $A(\lambda)$.
\end{lem}

\begin{proof}
Let $R(\lambda)=\colim F_k(\lambda)$. 
By construction, $A(\lambda)=\bigcup\limits_{\lambda'\in
T(\lambda)} H(\lambda')$. 
Since the $H(\lambda')$ are closed subspaces of $A(\lambda)$,
if suffices to show that the canonical map
 $p:R(\lambda)\to A(\lambda)$ is bijective. It is clearly surjective. So let $x\in
H(\lambda_1)\cap H(\lambda_2)
\subset A(\lambda)$. We need to show that $x\in H(\lambda_1)$ is related to $x\in H(\lambda_2)$
in the colimit $R(\lambda)$. Now $x=(c_1,\ldots,c_k)\in \mathcal{C}_n(1)^k$, and the
little cubes $c_1,\ldots,c_k$
satisfy ordering conditions specified by $\lambda_1$ and $\lambda_2$. Define
$\lambda_3\in\widehat{\mathcal{K}}_n(k)$
as follows:
the edge between $p$ and $q$ obtains no color or
orientation if the corresponding edges in $\lambda_1$ and $\lambda_2$ are not colored (note:
by definition of $T_k(\lambda)$ 
an edge in $\lambda_1$ is not colored iff the corresponding edge in 
$\lambda_2$ is not colored).
If both  are
colored, the corresponding edge in $\lambda_3$ obtains the color and
orientation of the edge with the smaller
color (if the colors agree, so do the orientations; this is forced by the
ordering conditions for $c_1,\ldots , c_k$). The ordering conditions for
$c_1,\ldots , c_k$ also imply that $\lambda_3$ does not have monochrome
cycles. By construction, $\lambda_3\in T_k(\lambda)$ and
$\lambda_3\le\lambda_1$ and
$\lambda_3\le \lambda_2$, and $x\in H(\lambda_3)$. Hence
 $x\in H(\lambda_1)$ and $x\in H(\lambda_2)$
represent the same
 element in the colimit.\\
This argument also proves the second statement.
\end{proof}

\begin{rema}\label{4_3a}
If $\lambda$ is the complete graph with no colors, then $T_k(\lambda) = 
\scK_n(k)$, and Lemma \ref{4_3} gives, in the
terminology of Berger \cite{Be}, a
``cellular decomposition'' of $\mathcal{C}_n(k)$ 
over $\mathcal{K}_n(k)$. In \cite{Be} Berger claimed that the same 
construction gives a cellular decomposition of $\mathcal{C}_n(k)$ over
$\mathcal{K}_n^B(k)$ and used this to show that $|\mathcal{K}_n^B|$ is
an $E_n$-operad. The following 
example illustrates that this construction does not give such a cellular
 decomposition over $\mathcal{K}_n^B$:

Let $(c_1,c_2,c_3)\in \mathcal{C}_3(3)$ be the configuration with $c_1=
[0,\frac{1}{2}]\times [\frac{2}{3},1]\times [0,\frac{1}{3}]$, $c_2=
[0,1]\times [\frac{1}{3},\frac{2}{3}]\times [\frac{1}{3},\frac{2}{3}]$, and
$c_3=[\frac{1}{2},1]\times [0,\frac{1}{3}]\times [\frac{2}{3},1]$. Over 
$\mathcal{K}_3^B$ this configuration lies in the interior of the cells
$C_\alpha$ and $C_\beta$ where 

$$
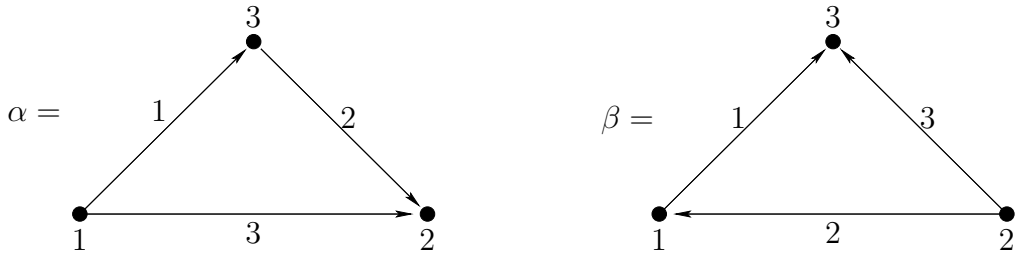
\begin{picture}(0,0)%
\includegraphics{berger1.pstex}%
\end{picture}%
\setlength{\unitlength}{2368sp}%
\begingroup\makeatletter\ifx\SetFigFont\undefined%
\gdef\SetFigFont#1#2#3#4#5{%
  \reset@font\fontsize{#1}{#2pt}%
  \fontfamily{#3}\fontseries{#4}\fontshape{#5}%
  \selectfont}%
\fi\endgroup%
\begin{picture}(10557,2589)(1051,-5608)
\put(2701,-4186){\makebox(0,0)[rb]{\smash{{\SetFigFont{12}{14.4}{\rmdefault}{\mddefault}{\updefault}{\color[rgb]{0,0,0}$1$}%
}}}}
\put(4501,-4261){\makebox(0,0)[lb]{\smash{{\SetFigFont{12}{14.4}{\rmdefault}{\mddefault}{\updefault}{\color[rgb]{0,0,0}$2$}%
}}}}
\put(3601,-5461){\makebox(0,0)[b]{\smash{{\SetFigFont{12}{14.4}{\rmdefault}{\mddefault}{\updefault}{\color[rgb]{0,0,0}$3$}%
}}}}
\put(8701,-4261){\makebox(0,0)[rb]{\smash{{\SetFigFont{12}{14.4}{\rmdefault}{\mddefault}{\updefault}{\color[rgb]{0,0,0}$1$}%
}}}}
\put(9601,-5461){\makebox(0,0)[b]{\smash{{\SetFigFont{12}{14.4}{\rmdefault}{\mddefault}{\updefault}{\color[rgb]{0,0,0}$2$}%
}}}}
\put(10501,-4261){\makebox(0,0)[lb]{\smash{{\SetFigFont{12}{14.4}{\rmdefault}{\mddefault}{\updefault}{\color[rgb]{0,0,0}$3$}%
}}}}
\put(3601,-3211){\makebox(0,0)[b]{\smash{{\SetFigFont{12}{14.4}{\rmdefault}{\mddefault}{\updefault}{\color[rgb]{0,0,0}$3$}%
}}}}
\put(7801,-5536){\makebox(0,0)[b]{\smash{{\SetFigFont{12}{14.4}{\rmdefault}{\mddefault}{\updefault}{\color[rgb]{0,0,0}$1$}%
}}}}
\put(11401,-5536){\makebox(0,0)[b]{\smash{{\SetFigFont{12}{14.4}{\rmdefault}{\mddefault}{\updefault}{\color[rgb]{0,0,0}$2$}%
}}}}
\put(9601,-3211){\makebox(0,0)[b]{\smash{{\SetFigFont{12}{14.4}{\rmdefault}{\mddefault}{\updefault}{\color[rgb]{0,0,0}$3$}%
}}}}
\put(1801,-5536){\makebox(0,0)[b]{\smash{{\SetFigFont{12}{14.4}{\rmdefault}{\mddefault}{\updefault}{\color[rgb]{0,0,0}$1$}%
}}}}
\put(5401,-5536){\makebox(0,0)[b]{\smash{{\SetFigFont{12}{14.4}{\rmdefault}{\mddefault}{\updefault}{\color[rgb]{0,0,0}$2$}%
}}}}
\put(7201,-4261){\makebox(0,0)[lb]{\smash{{\SetFigFont{12}{14.4}{\rmdefault}{\mddefault}{\updefault}{\color[rgb]{0,0,0}$\beta =$}%
}}}}
\put(1051,-4186){\makebox(0,0)[lb]{\smash{{\SetFigFont{12}{14.4}{\rmdefault}{\mddefault}{\updefault}{\color[rgb]{0,0,0}$\alpha =$}%
}}}}
\end{picture}%

$$

So the cells $C_\alpha$ and $C_\beta$ do not have disjoint interiors, which
violates Berger's notion of a cellular decomposition.
In contrast to $\mathcal{K}^B_3$, over $\mathcal{K}_3$,
 this configuration lies in $C_\gamma$ with

$$
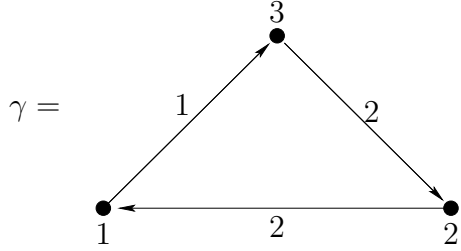
\begin{picture}(0,0)%
\includegraphics{berger2.pstex}%
\end{picture}%
\setlength{\unitlength}{2368sp}%
\begingroup\makeatletter\ifx\SetFigFont\undefined%
\gdef\SetFigFont#1#2#3#4#5{%
  \reset@font\fontsize{#1}{#2pt}%
  \fontfamily{#3}\fontseries{#4}\fontshape{#5}%
  \selectfont}%
\fi\endgroup%
\begin{picture}(4782,2589)(826,-5608)
\put(2701,-4186){\makebox(0,0)[rb]{\smash{{\SetFigFont{12}{14.4}{\rmdefault}{\mddefault}{\updefault}{\color[rgb]{0,0,0}$1$}%
}}}}
\put(4501,-4261){\makebox(0,0)[lb]{\smash{{\SetFigFont{12}{14.4}{\rmdefault}{\mddefault}{\updefault}{\color[rgb]{0,0,0}$2$}%
}}}}
\put(3601,-3211){\makebox(0,0)[b]{\smash{{\SetFigFont{12}{14.4}{\rmdefault}{\mddefault}{\updefault}{\color[rgb]{0,0,0}$3$}%
}}}}
\put(3601,-5461){\makebox(0,0)[b]{\smash{{\SetFigFont{12}{14.4}{\rmdefault}{\mddefault}{\updefault}{\color[rgb]{0,0,0}$2$}%
}}}}
\put(1801,-5536){\makebox(0,0)[b]{\smash{{\SetFigFont{12}{14.4}{\rmdefault}{\mddefault}{\updefault}{\color[rgb]{0,0,0}$1$}%
}}}}
\put(5401,-5536){\makebox(0,0)[b]{\smash{{\SetFigFont{12}{14.4}{\rmdefault}{\mddefault}{\updefault}{\color[rgb]{0,0,0}$2$}%
}}}}
\put(826,-4186){\makebox(0,0)[lb]{\smash{{\SetFigFont{12}{14.4}{\rmdefault}{\mddefault}{\updefault}{\color[rgb]{0,0,0}$\gamma =$}%
}}}}
\end{picture}%

$$

 which is in the
boundary of $C_\alpha$ and $C_\beta$.
For definitions and terminology consult \cite{Be}.

We want to point out that results from \cite{Be} and \cite{BFSV} imply that
$|\mathcal{K}_n^B|$ is an $E_n$-operad and that the inclusion $|\mathcal{K}_n^B|
\subset |\mathcal{K}_n|$ is
a $\Sigma$-equivalence.
\end{rema}

The colimit decompositions of the $A(\lambda)$ are functorial with respect to the 
colimit decomposition of $(\Ass\otimes\mathcal{C}_n)(k)$ of Lemma \ref{4_2} in the 
following sense: $T_k$ defines a functor
$$
T_k: \widehat{\mathcal{K}}_1(k)^{\op} \longrightarrow \mathcal{P}osets, \qquad
\lambda \mapsto T_k(\lambda)
$$ 
For a morphism $f:\lambda_1\to\lambda_2$ in
$\widehat{\mathcal{K}}_1(k)^{\op}$, i.e. $\lambda_2\le\lambda_1$, we have
a map of posets $T_k(f):T_k(\lambda_1)\longrightarrow T_k(\lambda_2)$ defined as follows:
let $S_i$ be the set of uncolored unoriented
edges of $\lambda_i$. Then $S_2\subset S_1$, and $E(k)\backslash
S_1$ is the set of uncolored edges of any $\mu\in T(\lambda_1)$. The map $T_k(f)$ sends $\mu$ to 
$\overline{\mu}$ obtained from $\mu$ by forgetting the orientations and colors of all edges in 
$E(k)\backslash S_2$. Moreover, $F_k(\lambda_1)(\mu)\subset F_k(\lambda_2)(\overline{\mu})$ because we
have less order conditions on the cubes in $F_k(\lambda_2)(\overline{\mu})$.
Hence the collection of functors $\{ F_k(\lambda);\ \lambda \in \widehat{\scK}_1(k)^{op}
\}$ is a $T_k$-indexed family of functors in the sense of \ref{A_1} below 
and we can combine the diagrams to a diagram
$$
F_k:\widehat{\mathcal{K}}_1(k)^{\op}\int T_k\longrightarrow \Top, \qquad (\lambda,\mu)
\mapsto F_k(\lambda)(\mu)
$$
where $\widehat{\mathcal{K}}_1(k)^{\op}\int T_k$ is the Grothendieck construction. Its
objects are pairs $(\lambda,\lambda')$  with $\lambda\in
\widehat{\mathcal{K}}_1(k)^{\op}$ and $\lambda'\in T_k(\lambda)$, 
and morphisms $(f,g):(\lambda_1,\lambda'_1)\to(\lambda_2,\lambda'_2)$
with $f:\lambda_1\to\lambda_2$ in $\widehat{\mathcal{K}}_1(k)^{\op}$, i.e.
$\lambda_2\le\lambda_1$ in
$\widehat{\mathcal{K}}_1(k)$,
 and $g:T_k(f)(\lambda'_1)=\overline{\lambda'_1}\to\lambda'_2$, i.e.
 $\overline{\lambda'_1}
\le\lambda'_2$, in
 $\widehat{\mathcal{K}}_n(k)$. 

\begin{lem}\label{4_a}
$(\Ass\otimes\mathcal{C}_n)(k)=\colim F_k$
\end{lem}

\begin{proof}
For $\lambda\in\widehat{\mathcal{K}}_1(k)^{\op}$ let
$$
i_\lambda: T_k(\lambda) \longrightarrow\widehat{\mathcal{K}}_1(k)^{\op}\int T_k
$$
denote the inclusion. Then
$$
\colim\nolimits_{\widehat{\mathcal{K}}_1(k)^{\op}\int T_k}
 F_k=\colim (\widehat{\mathcal{K}}_1(k)^{\op}\to\Top, \ \lambda\mapsto 
\colim\nolimits_{T_k(\lambda)}F_k\circ i_\lambda)
$$
Since $F_k\circ i_\lambda = F_k(\lambda)$ Lemmas \ref{4_2} and \ref{4_3} imply the statement.
\end{proof}

\begin{lem}\label{4_4}
There is an isomorphism of categories
$$
\xymatrix{
\varphi: \mathcal{K}_{n+1}(k) \ar[r] &\widehat{\mathcal{K}}_1(k)^{\op}\int T_k
}
$$
\end{lem}

\begin{proof}
The isomorphism is defined by sending $\lambda\in\mathcal{K}_{n+1}(k)$ to
$(\varphi_1(\lambda),\varphi_2
(\lambda))$ where $\varphi_1(\lambda)$ is obtained from
$\lambda$ by replacing the color $n+1$ by the color 1 and
by deleting colors and orientations of edges colored by any
$i\le n$, and $\varphi_2(\lambda)$ by
forgetting colors and orientations of all edges colored by $n+1$.
Observe that $\varphi_2(\lambda)$ is a complementary graph of $\varphi_1(\lambda)$.
\end{proof}

\begin{lem}\label{4_5}
For each $\lambda\in\mathcal{K}_{n+1}(k)$ the space $(F_k\circ\varphi) (\lambda)$ is contractible.
\end{lem}

\begin{proof}
$(F_k\circ\varphi)(\lambda)=\bigcup\{ H(\mu); \ \mu\in T_k(\varphi_1(\lambda)),\
\mu\le \varphi_2(\lambda) \}$.
Berger 
has proved the contractability of $\bigcup\{ H(\mu); \ \mu\in \mathcal{K}_n(k),\ \mu\le \lambda \}$
with $\lambda\in \mathcal{K}^B_n(k)$ \cite[Thm. 5.5]{Be}. The same argument applies to our situation.
\end{proof}

\begin{lem}\label{4_6}
The diagram $F_k\circ\varphi $
 is  Reedy-cofibrant, i.e. for all $\lambda\in\mathcal{K}_{n+1}(k)$
$$
\colim\nolimits_{\mu<\lambda}(F_k\circ\varphi)(\mu)\longrightarrow (F_k\circ\varphi)(\lambda)
$$
is a closed cofibration.
\end{lem}

\begin{proof}
By \ref{4_3}, $\colim_{\mu<\lambda}(F_k\circ\varphi)(\mu)=
\bigcup\limits_{\mu<\lambda} H
(\varphi_2(\mu))$, so that we have to show that
$$
\bigcup\limits_{\mu<\lambda} H(\varphi_2(\mu))\subset \bigcup\limits_{\mu\le\lambda} 
H(\varphi_2(\mu))
$$
is a closed cofibration. Let $\mu_1,\ldots,\mu_n$ be the predecessors of $\lambda$ and let 
$\lambda_i=\varphi_2
(\mu_i)$, $\lambda'=\varphi_2(\lambda)$. Since
$$
\xymatrix{
\bigcup^n_{i=1} H(\lambda_i)\cap H(\lambda') \ar[r]\ar[d] &
H(\lambda') \ar[d]
\\
\bigcup^n_{i=1} H(\lambda_i) \ar[r] & F_k(\lambda)
}
$$
is a pushout, it sufficies to show that
$$
\bigcup\nolimits^n_{i=1} H(\lambda_i)\cap H(\lambda')\longrightarrow H(\lambda')
$$
is a closed cofibration. By Lillig's union theorem \cite[Cor. 3]{Lillig} this holds if 
for any choice of
objects $\mu_1,\ldots,\mu_r\in\widehat{\mathcal{K}}_n(k)$ the map
$$
H(\mu_1)\cap\ldots\cap H(\mu_r)\longrightarrow H(\mu_r)
$$
is a closed cofibration. A little cube $c\in\mathcal{C}_n(1)$ is determined by its 
lowest vertex $x=(x_1
,\ldots, x_n)$ and its highest vertex $y=(y_1,\ldots,y_n)$. So $\mathcal{C}_n(1)^k\subset 
\mathbb{R}^{2n
k}$ is given by inequalities
$$
0\le x_{ij}< y_{ij}\le 1 \quad i=1,\ldots,k,\;j=1,\ldots,n.
$$
The subspace $H(\mu)\subset \mathcal{C}_n(1)^k$ consists of elements satisfying 
additional non-strict inequalities given by the ordering conditions.

Let $A\subset \mathbb{R}^{2nk}$ be the subspace given by all inqualities determining 
$H(\mu_1)\cap\ldots
\cap H(\mu_r)$ made non-strict, and $X$ the corresponding space obtained from $H(\mu_r)$. 
Then $A\subset X$
clearly is a closed cofibration. Define $\tau: X\to[0,1]$ to be the product of all
$(y_{ij}-x_{ij})$ for
 which we have strict inequalities in $H(\mu_1)\cap\ldots\cap H(\mu_r)$ 
(they are the same as the ones in
the list of $H(\mu_r)$), and let $V=\tau^{-1}(]0,1])$. Then by a result of Dold
 \cite[Satz 1]{Dold}
$$
V\cap A=H(\mu_1)\cap\ldots\cap H(\mu_r)\subset V=H(\mu_r)
$$
is a closed cofibration.
\end{proof}

\begin{coro}\label{4_7}
The canonical map $\hocolim (F_k\circ\varphi)\to 
\colim (F_k\circ\varphi)=(\Ass\otimes\mathcal{C}_n)(k)$ is a 
homotopy equivalence.
\end{coro}

For a proof see \cite[Prop. 6.9]{BFSV}. \hfill $\square$

\begin{lem}\label{4_8}
The operad $\Ass\otimes\mathcal{C}_n$ is $\Sigma$-free.
\end{lem}

\begin{proof}
By  \cite[Cor. 5.7]{BV3} the $\Sigma_k$-action on $(\Ass\otimes\mathcal{C}_n)(k)$ is free. 
Since each space
$(\Ass\otimes\mathcal{C}_n)(k)$ is Hausdorff and paracompact, the lemma follows.
\end{proof}

\begin{leer}\label{4_8a} The maps
$$
\xymatrix{
|\mathcal{K}_{n+1}(k)| &
\hocolim (F_k\circ\varphi) \ar[l]\ar[r] &
(\Ass\otimes\mathcal{C}_n)(k)
}
$$
assemble to $\Sigma$-equivalences of operads. 
\end{leer}
For each $k$ the maps are homotopy equivalences by Lemma \ref{4_5} and
Corollary \ref{4_7}. Both are equivariant homotopy equivalences since
$|\mathcal{K}_{n+1}(k)|$ and $(\Ass\otimes\mathcal{C}_n)(k)$ are free
$\Sigma_k$-spaces. It remains to prove that the collection of these
maps form maps of operads. For this it sufficies to show that

\begin{leer}\label{4_9}
$H(\varphi_2(\lambda))\circ (H(\varphi_2(\lambda_1)\times\ldots\times 
H(\varphi_2 )(\lambda_k))\subset
H(\varphi_2(\lambda\circ(\lambda_1\oplus\ldots\oplus\lambda_k))$ 
\end{leer}
for 
$\lambda\in\mathcal{K}_{n+1}(k)$ and $\lambda_i\in\mathcal{K}_{n+1}(l_i), \ i=1,\ldots, k$.
On the left side, composition is determined by the one in $\widehat{\mathcal{K}}_1\int
\mathcal{C}_n)$, on the right side we have composition in $\mathcal{K}_{n+1}$.

Condition (\ref{4_9}) is a consequence of the following properties of 
our order relations on $\mathcal{C}_n(1)$:
\begin{leer}\label{4_10}
\begin{enumerate}
\item[(i)] $c_1<_i c_2$ iff $c_3\circ c_1<_i c_3\circ c_2$ for all $c_3\in\mathcal{C}_n(1)$
\item[(ii)] $c_1<_i c_2\Rightarrow c_1\circ c_3<_ic_2\circ c_4$ for all $c_3,c_4\in\mathcal{C}_n(1)$
\end{enumerate}
\end{leer}
Finally we show
\begin{lem}\label{4_11}
$|\mathcal{K}_n|$ is an $E_n$-operad for each $n$.
\end{lem}
\begin{proof}
$\varphi (\mathcal{K}_n(k))\subset \widehat{\scK}_1(k)^{op}\int T_k$ 
exactly consists of all those pairs $(\varphi_1(\lambda),
\varphi_2(\lambda))$ for which $ \varphi_1(\lambda)$ does nor have any colors.
By our arguments of Remark \ref{4_3a}, Diagram \ref{4_8a} restricts to a
diagram of $\Sigma_k$-equivariant homotopy equivalences
$$
\xymatrix{
|\mathcal{K}_n(k)| &
\hocolim (F_k\circ(\varphi|\mathcal{K}_n(k))) \ar[l]\ar[r] &
\mathcal{C}_n(k)
}
$$
\end{proof}

To distinguish between the $F_k$ for the various $n$ we denote $F_k$ above by $F^{(n)}_k$ and
similarly for $\varphi$. 
We summarize:
\begin{leer}\label{4_12}
There is an explicit chain of $\Sigma$-equivalences
$$\scA ss\otimes \scC_n\leftarrow \hocolim (F^{(n)}\circ \varphi^{(n)})\rightarrow
|\scK_{n+1}|\leftarrow \hocolim (F^{(n+1)}\circ (\varphi^{(n+1)}|\scK_{n+1}))\rightarrow
\scC_{n+1}
$$
\end{leer}
Together with the Observation \ref{4_1} this completes the proof of Theorem C.

\section{Appendix: Iterated colimits and the Gro\-then\-dieck construction}

Our description of $(\Ass\otimes \sC_n)(k)$ as an iterated colimit is a
special case of a more general situation, which may be of
separate interest.

Let $F:\scA\to\scC at$ be any functor. Recall that the Grothendieck
construction $\scA\int F$ is the category whose objects are pairs
$(A,B)$ with $A\in{\rm obj}\scA$ and with $B\in{\rm obj}F(A)$. A morphism
$(A_1,B_1)\longrightarrow(A_2,B_2)$ is a pair $(\alpha,\beta)$ where
$\alpha:A_1\to A_2$ and $\beta:F(\alpha)(B_1)\to B_2$.

\begin{defi}\label{A_1}
 An $F$-\textit{indexed family of functors} into a category $\scC$ is
a collection of functors
$$\{G_A:F(A)\longrightarrow\scC: \ A\in{\rm obj}\scA \}$$
and natural transformations
$$\{\eta_\alpha: G_{A_1}\longrightarrow G_{A_2}F(\alpha);
\ \alpha:A_1\to A_2\ \in{\rm mor}\scA\}$$
satisfying $\eta_{id_A}=id_{G_A}$ for all $A\in{\rm obj}\scA$ and satisfying the
the following associativity conditions
$$\xymatrix{
G_{A_1}\ar[rrr]^(.4){\eta_{\alpha_2\alpha_1}}\ar[drr]^{\eta_{\alpha_1}}
&&&G_{A_3}F(\alpha_2\alpha_1)=G_{A_3}F(\alpha_2)F(\alpha_1)\\
&&G_{A_2}F(\alpha_1)\ar[ur]_{\eta_2F(\alpha_1)}
}$$
for any composable pair of morphisms
$$A_1\stackrel{\alpha_1}{\longrightarrow}A_2\stackrel{\alpha_2}{\longrightarrow}A_3$$
in $\scA$.
\end{defi}

An $F$-indexed family of functors determines a functor
$G\int F:\scA\int F\longrightarrow\scC$ given on objects by
$G\int F(A,B)=G_A(B)$ and on morphisms $(\alpha,\beta):(A_1,B_1)\longrightarrow(A_2,B_2)$ by
$$G_{A_1}(B_1)\stackrel{\eta_\alpha}{\longrightarrow}G_{A_2}F(\alpha)(B_1)
\stackrel{G_{A_2}(\beta)}{\longrightarrow}G_{A_2}(B_2).$$

Now suppose $\scC$ is a category with small colimits. Then the natural
tranformations $\eta_{\alpha}$ induce a functor $\scA\to\scC$, which takes
an object $A$ to $\colim_{F(A)} G_A$. We then have
\bigskip

\begin{prop}\label{A_2} $
\colim_{A\in{\rm obj}\scA}\left(\colim_{F(A)} G_A\right)
\cong\colim_{\scA\int F}G\int F$
\end{prop}
The proof is straight forward.


\begin{thebibliography}{99}
\bibitem{BFSV} C. Balteanu, Z. Fiedorowicz, R. Schw\"anzl, and R. M. Vogt,
Iterated monoidal categories,
Adv. in Math. 176 (2003), 277--349.

\bibitem{BM} M. Basterra, M.A. Mandell, The multiplication on $MU$ and $BP$,
to appear.

\bibitem{Be} C. Berger,
Combinatorial models for real configuration spaces and $E_n$ operads,
Contemp. Math. 202 (1997), 37--52.

\bibitem{BV1} J.M. Boardman, R.M. Vogt, 
Homotopy-everything $H$-spaces, 
Bull. Amer. Math.  Soc. 74 (1988), 1117-1122. 

\bibitem{BV2} J.M. Boardman, R.M. Vogt, 
Homotopy invariant structures on topological spaces,
Lecture Notes in Math. 347, Springer Verlag, Berlin 1973.

\bibitem{BV3} J.M. Boardman, R.M. Vogt, 
Tensor products of theories, application to infinite loop spaces,
J. Pure Appl. Algebra 14 (1979), 117-129.

\bibitem{Dold}
A. Dold, Die Homotopieerweiterungseigenschaft (= HEP) ist eine lokale
Eigenschaft, Inventiones math. 6 (1968), 185-189.

\bibitem{EKMM}
A.D. Elmendorf, I. Kriz, M.A. Mandell, J.P. May,
Rings, modules, and algebras in stable homotopy theory, 
Mathematical Surveys Monograms 47, Amer. Math. Soc., Providence, RI, 1996.

\bibitem{F}
Z. Fiedorowicz, Constructions of $E_n$ Operads, Proceedings of the Workshop
on Operads (1999), 34-55, University of Bielefeld.

\bibitem{Kon}
M. Kontsevich, 
Operads and motives in deformation quantification,
Lett. Math. Phys. 48 (1999), 35-72.

\bibitem{KS}
M. Kontsevich, Y. Soibelman, Deformations of algebras over operads and the
Deligne conjecture, Conference Moshe Flato 1999, vol. I, Math.
Phys. Stud. 21 (2000), 255-307.

\bibitem{Lillig}
J. Lillig, A union theorem for cofibrations, Arch. Math. 24 (1973), 410-415.

\bibitem{MSV}
J.E. McClure, R. Schw"anzl, R.M. Vogt,
$THH(R)\cong R\otimes S^1$ for $E_\infty$ ring spectra,
J. Pure Appl. Algebra 140 (1990), 23-32.

\bibitem{MS}
J.E. McClure, J.H. Smith, 
A solution of Deligne's Hochschild cohomology conjecture, Contemp.
Math. 293 (2002), 153-193.

\bibitem{Tam1}
D. Tamarkin, 
Another proof of M. Kontsevich formality theorem, 
Preprint (1998), math.QA/9803025.

\bibitem{Tam2}
D.E. Tamarkin, 
Formality of chain operad of little discs, Lett. Math. Phys. 66
(2003), 65-72.

\bibitem{Vogt}
R.M. Vogt,
Cofibrant operads and $E_\infty$-operads,
Topology and it Applications 133 (2003), 69-87.

\bibitem{Vor}
A.A. Voronov,
Homotopy Gerstenhaber algebras, Conference Moshe Flato 1999, vol. II, Math.
Phys. Stud. 22 (2000), 307-331.
\end{thebibliography}
\end{document}